\documentstyle[12pt,amscd]{amsart}

\newtheorem{Theorem}{Theorem}[section] 
\newtheorem{Lemma}[Theorem]{Lemma}
\newtheorem{Corollary}[Theorem]{Corollary}
\newtheorem{Proposition}[Theorem]{Proposition}

\newtheorem{Question}[Theorem]{Question}

\def\To{\longrightarrow}
\def\Proj{\operatorname{Proj}}

\def\rank{\operatorname{rank}}

\def\height{\operatorname{ht}}

\def\a{{\alpha}}
\def\b{{\beta}}
\def\l{{\lambda}}
\def\e{{\varepsilon}}
\def\NN{{\mathbb N}}
\def\ZZ{{\mathbb Z}}
\def\RR{{\mathbb R}}
\def\CC{{\mathbb C}}
\def\PP{{\mathbb P}}

\def\Q{{\bf Q}}

\def\m{{\mathfrak m}}
\def\aa{{\mathfrak a}}

\begin{document}

\title{Mixed multiplicities of ideals\\ versus mixed volumes of polytopes}
\author{Ngo Viet Trung and Jugal Verma }
\address{Institute of Mathematics,  Vi\^en To\'an Hoc, 
         18 Ho\`ang Qu\^oc Vi\^et, 10307 Hanoi, Vietnam}
\email{nvtrung@@math.ac.vn}
\address{ Department of Mathematics, Indian Institute of Technology Bombay, 
          Mumbai, India 400076}
\email{jkv@@math.iitb.ac.in}

\begin{abstract}
The main results of this paper interpret mixed volumes of lattice 
polytopes as mixed multiplicities of ideals and mixed multiplicities of ideals as 
Samuel's multiplicities. 
In particular, we can give a purely algebraic proof of Bernstein's theorem 
which asserts that the number of common zeros of a system of Laurent 
polynomial equations in the torus is bounded above by the mixed volume of 
their Newton polytopes. 

\end{abstract}

\maketitle

\thispagestyle{empty}
\section*{Introduction}

Let  us first recall the definition of mixed volumes.
Given two polytopes  $P, Q$ in $\RR^n$ (which need not to be 
different), their Minkowski sum is defined as the polytope
$$P + Q :=\{a + b \mid \ a \in P,\ b \in Q\}$$ 
The $n$-dimensional {\it mixed volume} of a collection of $n$ 
polytopes $Q_1,...,Q_n$ in 
$\RR^n$ is the value 
$$MV_n(Q_1, \ldots,Q_n) := 
\sum_{h=1}^s \sum_{1\le i_1 <...< i_h\le n} (-1)^{n-h}
  V_n(Q_{i_1}+\cdots+Q_{i_h}).$$
Here $V_n$ denotes the $n$-dimensional Euclidean volume.
Mixed volumes play an 
important role in convex geometry (see [BF], [Ew]) and elimination theory (see [GKZ], [CLO], [Stu]). \par

Our interest in mixed volumes arises from the following result of Bernstein 
[Be]  which relates the number of solutions of a system of polynomial 
equations to the mixed volume of their Newton polytopes 
(see also [Kh], [Ku]). \smallskip

\noindent{\bf Bernstein's Theorem.} 
Let $f_1,...,f_n$ be Laurent polynomials in 
$\CC[x_1^{\pm 1},...,x_n^{\pm 1}]$  with 
finitely many common zeros in the torus $(\CC^*)^n.$  Then the number of 
common 
zeros of $f_1,...,f_n$ in $(\CC^*)^n$ is bounded above by the mixed volume
$MV_n(Q_1,...,Q_n)$, where $Q_i$ denotes the Newton polytope of $f_i$. 
Moreover, this bound is attained  for a generic choice of coefficients in 
$f_1,..., f_n$.  
\smallskip

Bernstein's theorem is a generalization of the classical Bezout's theorem. 
It is a beautiful example of the interaction between  
algebra and combinatorics. However, the original proof in [Be] has more or 
less a combinatorial flavor. A geometric proof using intersection theory was 
given by Teissier [T3] (see also the expositions [Fu], [GKZ]). This paper 
grew out of our attempt to find an algebraic proof of Bernstein's theorem by using Samuel's multiplicity as it is usually done in a proof of Bezout's theorem. The relationship between toric varieties and multigraded rings used in the geometric proof suggests that 
mixed multiplicities of ideals may be the link between 
mixed volume of Newton polytopes of Laurent 
polynomials and the  number of their common zeros. 
To produce this link we encountered two problems 
which are of independent interest:

\begin{itemize}
\item Can one interpret the number of common zeros of Laurent polynomials 
in the torus as mixed multiplicity of ideals?
\item Does there exist any relationship between mixed 
multiplicities of ideals and mixed volume of polytopes?
\end{itemize}

We will solve these problems  and we will obtain thereby a  proof
for Bernstein's theorem which uses mixed multiplicities of ideals in a similar way as Samuel's multiplicity for Bezout's theorem. 
In fact, the number of common zeros of general polynomials in the torus 
counted with multiplicities and 
the mixed volume of their Newton polytopes can be interpreted as the same 
mixed multiplicity of ideals. \smallskip

Now we are going to give a brief introduction of mixed multiplicities. 
Let $J_1,...,J_n$ be a collection of ideals in a local ring $(A,\m)$ and 
$I$ an $\m$-primary ideal. Then the length function 
$$\ell(I^{u_0}J_1^{u_1}\cdots J_n^{u_n}/I^{u_0+1}J_1^{u_1}\cdots J_n^{u_n})$$ 
is a polynomial $P(u)$ for $u_0, u_1,...,u_n$ large enough [Ba], [R2], [Te1].
If we write this polynomial in the form
$$P(u) = \sum_{\a \in \NN^{n+1},\;\; |\a| = r} \frac{1}{\a!}e_\a u^\a + 
\{\text{terms of total degree $< r$}\},$$
where $r = \deg P(u)$ and $\a = (\a_0,\a_1,...,\a_n)$ of weight
$$
|\a| := \a_0 +\a_1 + \cdots + \a_n = r,\ \
\a!  := \a_0!\a_1!...\a_n!,\ \
u^\a  := u_0^{\a_0}u_1^{\a_1}...u_n^{\a_n},
$$
then the coefficients $e_\a$ are non-negative integers. One calls $e_\a$ the 
{\it mixed multiplicities} of the ideals $I,J_1,\ldots,J_n$ [Te1]. We will denote 
$e_\a$ by $e_\a(I|J_1,...,J_n)$. This notion can  also be defined for 
homogeneous ideals in a standard multi-graded algebra over a field.
Applications of mixed multiplicities can be found in 
[KaV], [Ro], [Te1], [Te2], [Tr2], [V1] and [V2].  \par

If the ideals $J_1,...,J_n$  are $\m$-primary ideals, one can interpret
$e_\a(I|J_1,...,J_n)$  as Samuel's multiplicity of general 
elements ([Te1], [R2], [Sw]). However, the techniques used in the $\m$-primary 
case are not applicable for non $\m$-primary ideals. For instance, mixed multiplicities of $\m$-primary ideals are always positive, whereas they may be zero in the general case. 
We will develop new techniques to prove the following general result which allows us to test the positivity of mixed multiplicities and  to compute them by means of Samuel's multiplicity.
\smallskip

\noindent{\bf Corollary \ref{main1}.}
Assume that the local ring $A$ has an infinite residue field. Let $Q$ be an ideal 
generated by $\a_i$ general elements in $J_i$, for $i = 1,...,n$, and 
$J := J_1 \cdots J_n$. Then $e_\a(I|J_1,...,J_n) > 0$ if and only 
if $\dim A/(Q:J^\infty) = \a_0+1$. In this case,
$$e_\a(I|J_1,...,J_n) = e(I,A/(Q:J^\infty)).$$

More generally, we can specify a class of concrete ideals $Q$ that can be used to compute $e_\a(I|J_1,...,J_n)$ (Theorem \ref{positivity}). Such a result was already obtained for two ideals in [Tr2]. The novelties here are the use of diagonal 
subalgebras and the introduction of superficial sequences for a set of ideals 
which provide us a simpler way to study mixed multiplicities. 
As consequences, we will show that the positivity of $e_\a(I|J_1,...,J_n)$ does not depend on the ideal 
$I$ and is rigid with respect to certain order of the indices $\a$ 
(Corollary \ref{rigid}). \par

There is already a close relationship between mixed multiplicities of multigraded rings
and mixed volumes.  Firstly, the multiplicity of a graded toric ring can 
be expressed in terms of the volume of a convex polytope (which is a 
consequence of Ehrhart's theory on the number of lattice points in 
convex polytopes). Secondly, mixed volume can be 
defined as a coefficient of the multivariate polynomial representing the 
volumes of linear combinations of the polytopes (Minkowski formula). 
Using these facts we find the following interpretation of mixed volumes as mixed multiplicities of ideals. \smallskip

\noindent {\bf Corollary \ref{main2}.}
Let $Q_1,...,Q_n$ be an arbitrary collection of lattice convex polytopes 
in $\RR^n$. Let $A = k[x_0,x_1,...,x_n]$ and $\m$ the maximal graded 
ideal of $A$. Let $M_i$ be any set of monomials of the same degree in 
$A$ such that $Q_i$ is the convex hull of the lattice points of their 
dehomogenized 
monomials in $k[x_1,...,x_n]$. Let $J_i$ be the ideal of $A$ generated 
by the monomials  of  $M_i$.  Then
$$MV_n(Q_1,...,Q_n) = e_{(0,1,...,1)}(\m|J_1,...,J_n).$$

This interpretation has interesting consequences. For instance, one can deduce properties of mixed volumes from those of mixed multiplicities. Conversely, properties of mixed volumes may predict unknown properties of mixed multiplicities. For instance, well-known inequalities for mixed 
volumes such as the Alexandroff-Fenchel inequality 
(see e.g. [Kh], [Te3]) lead
us to raise the question whether similar inequalities are valid for 
mixed multiplicities of ideals (Question \ref{Hodge}). To give an answer to 
this question turns out to be a challenging problem. \par

To prove Bernstein's theorem we first reformulate it for a system 
of homogeneous polynomial equations. In this case, the number of 
common zeros  of  general polynomials $f_1,\ldots,f_n$ can be seen as the Samuel's multiplicity of certain graded algebra. It turn out that this Samuel's multiplicity and the mixed volume of 
their Newton polytopes are the same mixed multiplicity $e_{(0,1,...,1)}(\m|J_1,...,J_n)$, where $J_1,...,J_n$ are the ideals generated by the supporting monomials of $f_1,\ldots,f_n$. By the 
principle of conservation of number, this implies the bound in Bernstein's 
theorem for any algebraically closed field (Theorem \ref{Bernstein}). \par

Finally, we would like to point out that computing  mixed volumes is a 
hard enumerative problem (see e.g. [EC], [HS1], [HS2]) and that the above 
relationships between mixed volumes, mixed multiplicities and Samuel 
multiplicity provide an alternative method for the computation of 
mixed volumes since many computer algebra programs can compute the 
Samuel multiplicity or the Hilbert polynomial of multigraded algebras. \par

This paper is organized as follows. Section 1 will deal with the 
characterization of mixed multiplicities as Samuel's multiplicities. 
In Section 2 we will interpret mixed volumes as mixed multiplicities. The algebraic proof of 
Bernstein's theorem  will be given in Section 3.\smallskip

\noindent
{\bf Acknowledgments.} The work on this paper began while the second author
visited the Institute of Mathematics, Hanoi in August 2002 under the
{\em China-India-Vietnam Network on Commutative Algebra and Algebraic Geometry} established by the Abdus Salam 
International Centre for Theoretical Physics (ICTP), Trieste, Italy. 
He thanks both institutions for their supports.  

\section{Mixed multiplicities of ideals}

We begin with some general observations on Hilbert polynomials of multigraded 
algebras.\par

Let $s$ be any non-negative integer. Let $R = \oplus_{u \in \NN^{s+1}}R_u$ 
be a  finitely generated standard $\NN^{s+1}$-graded algebra 
over an Artin local ring $R_0.$ We say $R$ is  {\it standard} if it is 
generated by homogeneous  elements of degrees $(0,..,1,..,0)$, where $1$ 
occurs only as the $i$th component, $i = 0,1,...,s$. The {\it Hilbert 
function} of $R$ is defined by $H_R(u) := \ell(R_u),$ 
where $\ell$ denotes the length.
If we view $u$ as a set of $s+1$ variables $u_0,...,u_s$, then there 
exists a polynomial $P_R(u)$ and integers $n_0,n_1,...,n_s$ such that 
$H_R(u) = P_R(u)$ for  $u_i \ge n_i$, $i = 0,1,...,s$ (abbr.~for $u \gg 0$) [Wa]. 
One calls $P_R(u)$ the {\it Hilbert polynomial}  of 
$R$.
If $P_R(u) \neq 0$, we write $P_R(u)$ in the form 
$$P_R(u) = \sum_{\a \in \NN^{s+1}, |\a| = r} \frac{1}{\a!}e_\a(R)u^\a + 
\{\text{terms of degree $< r$}\},$$
where $r = \deg P_R(u)$ and $\a = (\a_0,\a_1,...,\a_s)$ with 
$$
|\a| := \a_0 +\a_1 + \cdots + \a_s = r, \ \
\a!  := \a_0!\a_1!...\a_s!,\ \ \mbox{and} \ \
u^\a  := u_0^{\a_0}u_1^{\a_1}...u_s^{\a_s}.
$$
One calls the coefficients $e_\a(R)$ the {\it mixed multiplicities} of 
the multigraded algebra $R$. 
If $s=0$, i.e. $R$ is an $\NN$-graded algebra, then $R$ has only one 
mixed multiplicity. It is the usual multiplicity of $R$ and we will denote 
it by $e(R)$. \par

The mixed multiplicities of $R$ can be studied by means  of  certain 
$\NN$-graded subalgebras.
Let $\l = (\l_0,\l_1,...,\l_s)$ be any sequence of non-negative integers. 
Set $$R^\l := \bigoplus_{n \ge 0}R_{n\l}.$$
Then $R^\l$ is a finitely  generated $\NN$-graded algebra over $R_0$. 
One calls $R^\l$ the $\l$-{\it diagonal subalgebra} of $R$. This notion plays 
an important role in the study of embeddings of blowups of projective 
schemes  [CHTV]. \par

\begin{Lemma} \label{diagonal} 
Let $r = \deg P_{R}(u) \ge 0$ and  let all components of $\l$ be  positive. Then  
$\dim R^\l = r+1$ and
$$e(R^\l) = r!\displaystyle \sum_{\a \in \NN^{s+1},\;\; |\a| = r} 
\frac{1}{\a!}e_\a(R)\l^\a.$$
\end{Lemma} 

\begin{pf}
Since all components of $\l$ are positive, we have
$$
P_{R^\l}(n)  = H_{R^\l}(n) = H_R(n\l)
= \sum_{\a \in \NN^{s+1}, \;\; |\a| = r} \frac{1}{\a!}e_\a(R)\l^\a n^r + 
\{\text{terms of degree $< r$\}}
$$
for $n\gg 0$. This implies the conclusion because $\dim R^\l = \deg 
P_{R^\l}(n)+1$.
\end{pf}

Let $(A,\m)$ be a local ring (or a standard graded algebra over a 
field, where $\m$ is the maximal graded ideal). Let $I$ be an $\m$-primary 
ideal and $J_1,\ldots,J_s$ a sequence of ideals of $A$.
One can define the $\NN^{s+1}$-graded algebra
$$R(I|J_1,\ldots,J_s) := \bigoplus_{(u_0,u_1,...,u_s) \in 
\NN^{s+1}}I^{u_0}J_1^{u_1}...J_s^{u_s}/I^{u_0+1}J_1^{u_1}...J_s^{u_s}.$$
This algebra can be viewed as the associated graded ring of 
the Rees algebra $A[J_1t_1,...,J_st_s]$ with respect to the ideal 
generated by the elements of $I$. \par

For short, set $R = R(I|J_1,\ldots,J_s)$. Then $R$ is a standard 
$\NN^{s+1}$-graded algebra. Hence 
it has a Hilbert polynomial $P_R(u)$. For any $\a \in \NN^{s+1}$ with 
$|\a| = \deg P_R(u)$ we will set 
$$e_\a(I|J_1,\ldots,J_s) := e_\a(R).$$
The mixed multiplicities $e_\a(I|J_1,\ldots,J_s)$ were studied first 
for $\m$-primary ideals in [Ba], [R1] [R2], [Te1]
and then for arbitrary ideals in [KaMV], [KaV], [Tr2], [Vi]. \par

Throughout this section let 
\begin{align*}
J & := J_1...J_s,\\
d & := \dim A/(0:J^\infty),
\end{align*}
where for any ideal $Q \subset A$ we set $Q:J^\infty :=  \cup_{m \ge 
0}(Q:J^m)$. Moreover, for any finitely generated $A$-module $E$ we will 
denote by  $e(I,E)$ the Samuel multiplicity of $E$ with respect to $I$. 

\begin{Theorem} \label{dimension}
Let $R = R(I|J_1,\ldots,J_s)$. Assume that $d = \dim A/(0:J^\infty)\ge 1$. 
Then \par
{\rm (a)} $\deg P_R(u) = d - 1$,\par
{\rm (b)} $e_{(d-1,0,...,0)}(I|J_1,\ldots,J_s) = e(I,A/(0:J^\infty))$.
\end{Theorem}

\begin{pf} 
Let $I',J_1',...,J_s'$ be the sequence of ideals generated by 
$I,J_1,...,J_s$ in the quotient ring $A/(0:J^\infty)$ and put 
$R' = R(I'|J_1',...,J_s')$. Then
\begin{align*}
R_u' &  = (I^{u_0}J_1^{u_1}...J_s^{u_s} + 
(0:J^\infty))/(I^{u_0+1}J_1^{u_1}...J_s^{u_s} + (0:J^\infty))\\
&  = I^{u_0}J_1^{u_1}...J_s^{u_s}/(I^{u_0+1}J_1^{u_1}...J_s^{u_s}
+ I^{u_0}J_1^{u_1}...J_s^{u_s}\cap (0:J^\infty)).
\end{align*}

Since $I^{u_0}J_1^{u_1}...J_s^{u_s}\cap (0:J^\infty) = 0$ for $u\gg 0$, 
we get
$R_u = R_u'$ for $u \gg 0$. Hence 
$$P_R(u) = P_{R'}(u).$$\par

So we may replace $A$ by $A/(0:J^\infty)$. 
If we do so, we may assume that $0:J^\infty = 0$ and $d = 
\dim A \ge 1$. Then $\height J \ge 1$.
For $\l = (1,...,1)$ we have 
$$R^\l = \oplus_{n\ge 0}I^{n}J^n/I^{n+1}J^n \cong A[IJt]/(I),$$
where $A[IJt]$ is the Rees algebra of the ideal $IJ$. 
Since $\height(IJ) \ge 1$, we have $\dim A[IJt] = d+1$ [Va,Corollary 1.6].
Hence $\dim R^\l \le d$. By Lemma \ref{diagonal} (a), this implies $\deg 
P_R(u) \le d-1$. \par

On the other hand, $\dim A/J^m < d$  for any $m \ge 1$. 
Therefore,
$$
e(I,A) = e(I,J^m) = \lim_{n \to \infty} 
\frac{\ell(I^nJ^m/I^{n+1}J^m)}{ n^{d-1}/(d-1)! } = \lim_{n\to \infty} 
\frac{P_R(n,m,...,m)}{n^{d-1}/(d-1)!}
$$
for $m \gg 0$. Since $e(I,A) > 0$, this implies $\dim P_R(u) \ge d-1$. 
So we can conclude that
$\deg P_R(u) = d-1$ and that
$e_{(d-1,0,...,0)}(R) =  e(I,A).$
\end{pf}

The computation of  mixed multiplicities can be passed to the case of 
$e_{(d-1,0,...,0)}(R)$. For this we shall need the following notation.\par

Given a standard $\ZZ^{s+1}$-graded algebra $S$, we will denote by 
$S_+$ the ideal of 
$S$ generated by the homogeneous elements of degrees with positive components. 
A sequence of homogeneous elements $z_1,...,z_m$ in $S$ is 
called {\it filter-regular} if
$$[(z_1,...,z_{i-1}):z_i]_u = (z_1,...,z_{i-1})_u$$
for $u \gg 0$, $i = 1,...,m$. It is easy to see that this is equivalent to 
the condition
$z_i \not\in P$ for any associated prime 
$P \not\supseteq S_+$ of $S/(z_1,...,z_{i-1})$. 

\noindent{\bf Remark.} {\rm Filter-regular sequences have their origin in the 
theory of Buchsbaum rings [SV, Appendix]. It can be shown that if $S$ is a 
standard graded algebra over a field, then $\Proj(S)$ is an equidimensional 
Cohen-Macaulay scheme if and only if every homogeneous system of parameters 
of $S$ is filter-regular.}
\smallskip

We will work now in the $\ZZ^{s+1}$-graded algebra
$$S := \bigoplus_{u\in\ZZ^{s+1}}I^{u_0}J_1^{u_1}...J_s^{u_s}
/I^{u_0+1}J_1^{u_1+1}...J_s^{u_s+1}.$$
Let $\e_1,...,\e_m$  be any non-decreasing sequence of indices with 
$1 \le \e_i \le s$. Let $x_1,...,x_m$ be a sequence of elements of $A$ 
with $x_i \in J_{\e_i}$, $i = 1,...,m$. We denote by $x_i^*$  
the residue class of $x_i$ 
in $J_{\e_i}/IJ_1...J_{\e_i-1}J_{\e_i}^2J_{\e_i+1}...J_s$.
We call $x_1,...,x_m$ an $(\e_1,...,\e_m)$-{\it superficial sequence} for the 
ideals $J_1,...,J_s$ (with respect to $I$) if $x_1^*,...,x_m^*$ is a 
filter-regular sequence in $S$. \par

The above notion can be considered as a generalization of  the classical notion of a 
superficial element of an ideal, which plays an important role in the theory 
of multiplicity. Recall that an element $x$ is called superficial with respect to an 
ideal $\aa$ if there is an integer $c$ such that 
$$(\aa^n:x) \cap \aa^c = \aa^{n-1}$$ for 
$n \gg 0$.  A sequence of elements $x_1,...,x_m \in \aa$ is called a 
superficial sequence of $\aa$ if the residue class of $x_i$ in 
$A/(x_1,...,x_{i-1})$ is a superficial element of the ideal 
$\aa/(x_1,...,x_{i-1})$, $i = 1,...,m$. It is known that this is equivalent to 
the condition that the initial forms of $x_1,...,x_m$ in $\aa/\aa^2$ 
form a filter-regular sequence in the associated graded ring 
$\oplus_{n\ge 0}\aa^n/\aa^{n+1}$ (see e.g. [Tr1, Lemma 6.2]).
\par

We may use superficial sequences to reduce the dimension of the base ring.

\begin{Lemma} \label{reduction}
Let $Q$ be an ideal of $A$ generated by an 
$(\e_1,...,\e_m)$-superficial sequence of 
$J_1,...,J_s$. Let $\bar I,\bar J_1,...,\bar J_s$ be the sequence of 
ideals generated  by $I,J_1,...,J_s$ in the quotient ring $A/Q$ and put 
$\bar R = R(\bar I|\bar J_1,...,\bar J_s)$. Let
$\a_j$ be the number of the indices $i$ such that $\e_i = j$, $j = 1,...,s$. 
Let  $\Delta^{(0,\a_1,...,\a_s)}P_R(u)$ denote the 
$(0,\a_1,...,\a_s)$-difference of the polynomial $P_R(u)$.
Then 
$$P_{\bar R}(u) = \Delta^{(0,\a_1,...,\a_s)}P_R(u),$$
\end{Lemma}

\begin{pf}
If $m = 1$, we may assume that $(\a_1,...,\a_s) = (1,0,...,0)$. 
Then $Q = (x)$, where $x \in J_1$ such that $(0:x^*)_u = 0$ for $u \gg 0$. 
This means
\begin{align}
(I^{u_0+1}J_1^{u_1+2}J_2^{u_2+1}...J_s^{u_s+1}:x) 
\cap I^{u_0}J_1^{u_1}...J_s^{u_s}  & = I^{u_0+1}J_1^{u_1+1}...J_s^{u_s+1} .
\end{align}
As a consequence we get
$$
(I^{u_0+1}J_1^{u_1+2}J_2^{u_2+1}...J_s^{u_s+1}:x) 
\cap I^{u_0}J_1^{u_1+1}...J_s^{u_s+1} = I^{u_0+1}J_1^{u_1+1}...J_s^{u_s+1}
$$
for $u \gg 0$. Consider $R$ as a quotient ring of $S$. The above 
formula shows that $(0_R:x^*)_u = 0$ for $u \gg 0$. 
Hence $P_{R/(0_R:x^*)}(u) = P_R(u).$
Now, from the exact sequence
$$0 \To R/(0 :_R x^*) \overset {x^*} {\To} R \To R/(x^*) \To 0$$
we can deduce that 
\begin{align} 
P_{R/(x^*)}(u) & = \Delta^{(0,1,0,...,0)}P_R(u).
\end{align}\par

On the other hand, (1) implies
$$I^{u_0+1}J_1^{u_1+2}...J_s^{u_s+1} 
\cap xI^{v_0}J_1^{v_1}...J_s^{v_s} = xI^{u_0+1}J_1^{u_1+1}...J_s^{u_s+1}$$
for $(u_0,u_1,...,u_s) \gg 0$ and $v_i = u_i, u_i+1$, $i = 0,1,...,s$. 
Using this formula we can easily show that
$$I^{u_0}J_1^{u_1}...J_s^{u_s} \cap xI^{v_0}J_1^{v_1}...J_s^{v_s} =
xI^{\max\{u_0,v_0\}}J_1^{\max\{u_1-1,v_1\}}...J_s^{\max\{u_s,v_s\}}$$
for $u \gg 0$, $v \gg 0$. 
By Artin-Rees  lemma, there exists 
$(c_0,c_1,...,c_s) \in \NN^{s+1}$ with $c_1 > 0$ such that
$$I^{u_0}J_1^{u_1}...J_s^{u_s} \cap (x) 
\subseteq xI^{u_0-c_0}J_1^{u_1-c_1}...J_s^{u_s-c_s}$$
for $u_i \ge c_i$, $i = 0,1,...,s$. Therefore,
\begin{align*} 
I^{u_0}J_1^{u_1}...J_s^{u_s} \cap (x) & = I^{u_0}J_1^{u_1}...J_s^{u_s} 
\cap xI^{u_0-c_0}J_1^{u_1-c_1}...J_s^{u_s-c_s}\\
& = xI^{u_0}J_1^{u_1-1}J_2^{u_2}...J_s^{u_s}
\end{align*}
for $u \gg 0$.  This implies
\begin{align*}
\bar R_u & = 
(I^{u_0}J_1^{u_1}...J_s^{u_s},x)/(I^{u_0+1}J_1^{u_1}...J_s^{u_s},x)\\
& = I^{u_0}J_1^{u_1}...J_s^{u_s}/(I^{u_0+1}J_1^{u_1}...J_s^{u_s} + 
I^{u_0}J_1^{u_1}...J_s^{u_s}\cap (x)),\\
& = 
I^{u_0}J_1^{u_1}...J_s^{u_s}/(I^{u_0+1}J_1^{u_1}\ldots 
J_s^{u_s} + xI^{u_0}J_1^{u_1-1}J_2^{u_2}...J_s^{u_s})\\
& = (R/(x^*))_u.
\end{align*}
Thus, $P_{\bar R}(u) = P_{R/(x^*)}(u).$
Combining this with (2) we get 
$P_{\bar R}(u) = \Delta^{(0,..,1,..,0)}P_R(u)$ 
which proves the case $m =1$.\par

If $m  > 1$, we may assume that $\a_1 > 0$. Then $x_1 \in J_1$. 
Let $I^*,J_1^*,...,J_s^*$ denote the sequence of the ideals 
generated by $I,J_1,...,J_s$ in the quotient ring $A/(x_1)$. 
Put $R^* = R(I^*|J_1^*,...,J_s^*)$. As shown above, we have
\begin{align}
P_{R^*}(u) & = \Delta^{(0,1,0,...,0)}P_R(u).
\end{align}
Let 
$S^* := \oplus_{u\in\ZZ^{s+1}}(I^*)^{u_0}(J_1^*)^{u_1} \ldots 
(J_s^*)^{u_s}/(I^*)^{u_0+1}(J_1^*)^{u_1+1}...(J_s^*)^{u_s+1}.$
For $u \gg 0$ we have
\begin{align*}
[S/(x_1^*)]_u & =  I^{u_0}J_1^{u_1}\ldots 
J_s^{u_s}/(I^{u_0+1}J_1^{u_1+1}\ldots J_s^{u_s+1}+x_1I^{u_0}J_1^{u_1-1}
J_2^{u_2}...J_s^{u_s})\\
& =  I^{u_0}J_1^{u_1}\ldots J_s^{u_s}/(I^{u_0+1}J_1^{u_1+1} \ldots 
J_s^{u_s+1}+(x_1) \cap I^{u_0}J_1^{u_1} \ldots J_s^{u_s})\\ 
& =  (I^{u_0}J_1^{u_1} \ldots J_s^{u_s},x_1)/(I^{u_0+1}J_1^{u_1+1} \ldots 
J_s^{u_s+1},x_1)\\
& = S^*_u,
\end{align*}
Since 
$[(x_1^*,...,x_{i-1}^*):x_i]_u = (x_1^*,...,x_{i-1}^*)_u$ 
for $u \gg 0$, $i = 2,...,m$, we also have
$$[(x_2^*,...,x_{i-1}^*)S^*:x_i]_u = (x_2^*,...,x_{i-1}^*)S^*_u$$
for $u \gg 0$, $i = 2,...,m$. Therefore, $x_2^*,...,x_m^*$  is an  
$(\e_1-1,\e_2,...,\e_s)$-superficial sequence  of the ideals 
$I^*,J_1^*,...,J_s^*$. Now, we may use induction on $m$ to assume that
$$P_{\bar R}(u) = \Delta^{(0,\a_1-1,\a_2,...,\a_s)}P_{R^*}(u).$$
Combining this with (3) we get 
$P_{\bar R}(u) = \Delta^{(0,\a_1,\a_2,...,\a_s)}P_R(u)$. 
\end{pf}

Using Theorem \ref{dimension} and Lemma \ref{reduction} we obtain the 
following criterion for the positivity of mixed multiplicities. 

\begin{Theorem} \label{positivity}
Let $\a = (\a_0,\a_1,...,\a_s)$ be any sequence of non-negative integers 
with $|\a| = d-1$.  Let $Q$ be any ideal generated by an 
$(\a_1,...,\a_s)$-superficial sequence of the ideals 
$I,J_1,...,J_s$.  Then $e_\a(I|J_1,...,J_s) > 0$ 
if and only if   $\dim A/(Q:J^\infty) = \a_0+1.$ In this case, 
$$e_\a(I|J_1,...,J_s) = e(I,A/(Q:J^\infty)).$$
\end{Theorem}

\begin{pf}
If $\a = (d-1,0,...,0)$, the conclusion follows from 
Theorem \ref{dimension}.
If $\a \neq (d-1,0,...,0)$, then $d \ge 2$. 
Let $\bar R, \bar I, \bar J_1,...,\bar J_s$ be as in 
Lemma \ref{reduction}. Then
$\deg P_{\bar R}(u) \le d-1-m = \a_0$ where $m=\a_1+ \cdots + \a_s.$ 
Write $$P_{\bar R}(u) = \sum_{\b \in \NN^{s+1}, |\b| 
= \a_0}\frac{e_\b(\bar I|\bar J_1,...,\bar J_s)}{\b!}u^\b + 
\{\text{terms of degree $< \a_0$}\}.$$
Then
$$e_{(\a_0,\a_1,...,\a_s)}(I|J_1,...,J_s) = 
e_{(\a_0,0,...,0)}(\bar I|\bar J_1,...,\bar J_s).$$\par
If $e_\a(I|J_1,...,J_s) > 0$ then 
$e_{(\a_0,0,...,0)}(\bar I|\bar J_1,...,\bar J_s) > 0$.  
Therefore, $\deg P_{\bar R}(u) = \a_0$. By Theorem \ref{dimension}(a),
this implies $\dim A/(Q:J^\infty) = \a_0+1$.\par
Conversely, if  $\dim A/(Q:J^\infty) = \a_0+1$ and if we put 
$\bar J = \bar J_1...\bar J_s$, then
$$e_{(\a_0,0,...,0)}(\bar I|\bar J_1,...,\bar J_s) 
= e(\bar I,\bar A/(0:\bar J^\infty)) = e(I,A/(Q:J^\infty))$$ 
by Theorem \ref{dimension}(b).
Since the Samuel multiplicity is always positive, 
this implies $e_{(\a_0,0,...,0)}(\bar I|\bar J_1,...,\bar J_s) > 0$.
So we can conclude that 
$e_\a(I|J_1,...,J_s) > 0$ if and only if $\dim A/(Q:J^\infty) = \a_0+1$.
\end{pf}

Let $k$ be the residue field of $A$. Using the prime avoidance 
characterization of a superficial element we can easily see that 
superficial sequences exist if $k$ is infinite. In fact, general 
elements of $J_1,...,J_s$ always form a superficial sequence. Recall that 
a property holds for a {\it general element} $x$ of an ideal 
$Q = (y_1,...,y_r)$ if there exists a non-empty Zariski-open 
subset $U \subseteq k^r$ such that whenever $x = \sum_{j=1}^m c_jx_j$ and 
the image of $(c_1,...,c_m)$ in $k^m$ belongs to $U$, then the  
property holds for $x$.

\begin{Lemma} \label{general}
Assume that $k$ is infinite. Any sequence which consists  
of $\a_1$ general elements in $J_1$, ... , $\a_s$ elements in $J_s$ 
forms an $(\a_1,...,\a_s)$-superficial sequence  for the ideals $J_1,...,J_s$.
\end{Lemma}

\begin{pf}
Let $x_1,...,x_m$ be a sequence of such general elements, 
$m = \a_1 + \cdots + \a_s$. Assume that $x_i \in J_{\e_i}$. 
Since $x_i$ is a general element of $J_{\e_i}$, 
we have $x_i^* \not\in P$ for any associated prime $P$ of 
$(x_1^*,...,x_{i-1}^*)$ with 
$P \not\supseteq J_{\e_i}/IJ_1...J_{\e_i-1}J_{\e_i}^2J_{\e_i+1}...J_s$. 
Since $S_+$ is contained in the ideal generated by the elements of 
$J_{\e_i}/IJ_1...J_{\e_i-1}J_{\e_i}^2J_{\e_i+1}...J_s$, this 
implies $x_i^* \not\in P$ for any associated prime $P$ of 
$(x_1^*,...,x_{i-1}^*)$ with $P \not\supseteq S_+$. 
Hence $x_1^*,...,x_m^*$ is a filter-regular sequence in $S$. 
\end{pf}

\begin{Corollary} \label{main1}
Assume that the local ring $A$ has infinite residue field. Let 
$Q$ be an ideal generated by $\a_1$ general elements in $J_1$, ... , $\a_n$ elements in $J_n$. Then $e_\a(I|J_1,...,J_n) > 0$ if and only if 
$\dim A/(Q:J^\infty) = \a_0+1$. In this case,
$$e_\a(I|J_1,...,J_n) = e(I,A/(Q:J^\infty)).$$
\end{Corollary}

Now we shall see that the characterization of mixed multiplicities of 
$\m$-primary ideals given in [Te1] is a special case of Corollary 
\ref{main1}.

\begin{Corollary} \label{Teissier} 
{\rm [Te1, Ch.~0, Proposition 2.1]}
Assume that the local ring $A$ has infinite residue field. 
Let $I,J_1,\ldots,J_s$ be $\m$-primary ideals. Let 
$\a = (\a_0,\a_1,...,\a_s)$ be any 
sequence of non-negative integers with $|\a| = \dim A-1$. Let 
$P$ be an ideal of $A$ generated by $\a_0+1$ general elements in 
$I$, $\a_1$ general elements in $J_1$, ... , $\a_s$ elements in $J_s$. Then
$$e_\a(I|J_1,...,J_s) = e(P,A).$$
\end{Corollary}

\begin{pf} 
Let $Q$ be the subideal of $P$ generated by $\a_1$ general elements in $J_1$, ... , $\a_s$ elements in $J_s$. By Lemma \ref{general}, 
these elements form a superficial sequence of the ideals $J_1,...,J_s$. 
Since $J_1,...,J_s$ are $\m$-primary ideals,  $Q$ is generated by a 
subsystem of parameters of $A$ and $J$ is an $\m$-primary ideal. 
Therefore,
$$\dim A/(Q:J^\infty) = \dim A/Q = \dim A-(\a_1+\cdots+\a_s) = \a_0+1.$$
By Theorem \ref{positivity}, and the above equation, we get
$$e_\a(I|J_1,...,J_s)  = e(I,A/(Q:J^\infty)) = e(I,A/Q).$$
But $e(I,A/Q) = e(P,A/Q)$ because $P$ generates a minimal reduction of 
$I$ in $A/Q$. So we can conclude that 
$$e_\a(I|J_1,...,J_s) = e(P,A/Q) = e(P,A).$$
\end{pf}

Using Corollary \ref{main1} we obtain interesting properties of 
mixed multiplicities.

\begin{Corollary} \label{rigid} 
Let $\a = (\a_0,\a_1,...,\a_s)$ be any sequence of non-negative integers 
with $|\a| = d -1$. Assume that $e_\a(I|J_1,\ldots,J_n) > 0$. 
Then\par
{\rm (a) } $e_\a(I'|J_1,\ldots,J_n) > 0$ for any $\m$-primary ideal 
$I'$,\par
{\rm (b) } $e_\b(I|J_1,\ldots,J_n) > 0$ for all 
$\b = (\b_0,\ldots,\b_n)$ with $|\b| = d-1$ 
and $\b_i \le \a_i$, $i = 1,\ldots,n$.
\end{Corollary}

\begin{pf}
Without loss of generality, we may assume that the residue field of $A$ 
is infinite.
Let $Q$ be an ideal generated by $\a_1$ general elements in $J_1$, ... , $\a_n$ elements in $J_n$. \par
(a) By Corollary \ref{main1}, the assumption implies 
$\dim A/(Q:J^\infty) = \a_0+1$. Since this condition does not depend on 
$I$, we also have $e_\a(I'|J_1,...,J_n) > 0$. \par
(b) Let $Q'$ denote the subideal of $Q$ generated by $\b_i$ 
general elements in $J_i$, $i = 1,...,n$. 
Put $A^* = A/Q'$, $I^* = IA^*$ and $J_i = J_iA^*$. 
Let $R^* = R(I^*|J_1^*,\ldots,J_n^*)$. By Lemma \ref{reduction} we have
$$P_{R^*}(u) = \Delta^{(0,\b_1,\ldots,\b_n)}P_R(u).$$
From this it follows that
$$e_{(\a_0,\a_1-\b_1,\ldots,\a_n-\b_n)}(I^*|J_1^*,\ldots,J_n^*)  
= e_\a(I|J_1,\ldots,J_n) > 0.$$
Hence $\deg P_{R^*}(u) = (d-1) - (\b_1 + \cdots +\b_n)$.
By Theorem \ref{dimension}(a), this implies
$$\dim A/(Q':J^\infty) = \deg P_{R^*}(u) + 1 = \b_0 + 1.$$
Therefore, $e_\b(I|J_1,\ldots,J_n) > 0$ by Corollary \ref{main1}.
\end{pf}

\section{Mixed volumes and toric rings}

The aim of this section is to interpret mixed volumes as 
mixed multiplicities.\par

Usually, mixed volume is defined for a collection of $n$ 
convex polytopes in $\RR^n$ (see e.g. [CLO]). But it is obvious that it may be also 
defined for any collection of convex polytopes in $\RR^n$ as follows. 
Let $Q_1,...,Q_r$ be convex polytopes in $\RR^n$ with 
$\dim (Q_1+\cdots+Q_r) \le r$. We call the value
$$ 
MV_r(Q_1, \ldots,Q_r) := 
\sum_{h=1}^r \sum_{1\le i_1 <...< i_h\le r} (-1)^{r-h}
  V_r(Q_{i_1}+\cdots+Q_{i_h}).
$$
the {\it mixed volume} of $Q_1,...,Q_r$. Here 
$V_r$ denotes the $r$-dimensional Euclidean  volume.
\par

Let $\Q = (Q_1,...,Q_s)$ be a sequence of convex polytopes in $\RR^n$. 
Let $\l = (\l_1,...,\l_s)$ be any sequence of non-negative integers.
We denote by $\l\Q$ the Minkowski sum $\l_1Q_1 + \cdots + \l_sQ_s$ 
and by
$\Q_\l$ the multiset of $\l_1$ polytopes $Q_1$,...,$\l_s$ polytopes 
$Q_s$.
Minkowski showed that the volume of the polytope $\l\Q$ is a 
homogeneous polynomial in $\l$ whose coefficients are mixed volumes up to 
constants  (see e.g. [CLO, Ch. 7, Proposition 4.9]).

\begin{Proposition} \label{Minkowski} {\rm (Minkowski formula)}
Let $r = \dim(Q_1+\cdots+Q_s)$. Then
$$V_r(\l\Q) = \sum_{\a \in \NN^s, |\a| = 
r}\frac{1}{\a!}MV_r(\Q_\a)\l^\a.$$
\end{Proposition}

We will use Minkowski formula to establish the relationship between
mixed volumes and mixed multiplicities. For that we need to work  
with graded toric rings.\par

Let $A = k[x_1,...,x_n]$ be a polynomial ring over a field $k$. Let
$M$ be a finite set of monomials in $A$. The subalgebra 
$k[M]$ of $A$ generated by the monomials of $M$ is called 
the {\it toric ring} (or affine semigroup ring) of $M$.
We associate with every monomial 
$x_1^{a_1}...x_n^{a_n} \in A$  the lattice point 
$a = (a_1,...,a_n) \in \NN^n$.
Many ring-theoretic properties of $k[M]$ can be described by means of 
the lattice points of $M$ (see e.g. [BH, Section 6] or [St, Chap. I]). 
For instance,  
$$\dim k[M] = \rank \ZZ(M),$$
where $\ZZ(M)$ denotes the subgroup of $\ZZ^n$ generated by the 
lattice points of $M$. \par

Assume furthermore that the lattice points of $M$ lie on an affine 
hyperlane of $\RR^n$.
This is for example the case when $M$  consists of monomials of the 
same degree. Then $k[M]$ has a natural $\NN$-graded structure. The 
multiplicity $e(k[M])$ can be expressed in terms of the lattice points of 
$M$ as follows.
\par

Let $Q_M$ denote the convex hull of the lattice points of $M$ in $\RR^n$.
Then $Q_M$ is a convex polytope with 
$$\dim Q_M = \rank \ZZ(M)-1.$$ 

\begin{Proposition}\label{volume} 
Let $r = \rank \ZZ(M)-1$. Let $E$ be any subset of $M$ such that 
its lattice points 
form a basis of $\ZZ(M)$. Then
$$e(k[M]) = \frac{V_r(Q_M)}{V_r(Q_E)}.$$
\end{Proposition}

This multiplicity formula is a consequence of Ehrhart's theory for the  
number of lattice points in  lattice polytopes 
(see  e.g. [BH, Theorem 6.3.12] or [St, Chap. I, Theorem 10.3]).  
The number $ V_r(Q_M)/V_r(Q_E)$ is often called the 
{\it normalized volume} of the polytope $Q_M$ with respect to the 
lattice $\ZZ(M)$.  \par

In the following we will be concerned with products of finite sets of 
monomials, which is the counterpart of Minkowski sums of convex 
polytopes.  \par

Let $M_1,...,M_s$ be sets of monomials in $A$ such that each $M_i$ 
consists of monomials of the same degree. For any sequence 
$\l = (\l_1,...,\l_s)$ of positive integers we denote by 
$M^\l$ the set of all products of $\l_1$ monomials of $M_1$,...,
$\l_s$ monomials  of $M_s$. 
Using the above propositions we can express the multiplicity of the 
toric ring $k[M^\l]$ in terms of mixed volumes.

\begin{Corollary} \label{product}
Let $r = \rank \ZZ(M^{(1,...,1)})-1$. Let $E$ be any subset of 
$M^{(1,...,1)}$ such that its lattice points form a basis of 
$\ZZ(M^{(1,...,1)})$. Let $\Q$ be the 
sequence of polytopes 
$Q_{M_1},...,Q_{M_s}$. Then
$$e(k[M^\l]) = \frac{1}{V_r(Q_E)}\sum_{\a \in \NN^s, |\a| = 
r}\frac{1}{\a!}MV_r(\Q_\a)\l^\a$$.
\end{Corollary}

\begin{pf}
Every lattice vector of $M^\l$ is a sum of $\l_1$ 
lattice points of $M_1$,...,$\l_s$ lattice points of $M_s$. Therefore,
$$Q_{M^\l} = \l_1Q_1 + \dots + \l_sQ_s = \l\Q.$$
Since $\ZZ(M^\l) = \ZZ(M^{(1,...,1)})$, we have $\rank \ZZ(M^\l) = r+1$. Using Proposition \ref{volume} we obtain 
$$e(k[M^\l]) =  \frac{V_r(\l\Q)}{V_r(Q_E)}.$$
Hence the conclusion follows from Proposition \ref{Minkowski}.
\end{pf}

This formula for the multiplicity of the toric rings $k[M^\l]$ 
resembles the formula for the multiplicity of diagonal subalgebras in 
Lemma \ref{diagonal}. Therefore, if we can find a standard multigraded 
algebra such that the toric rings $k[M^\l]$ are its diagonal subalgebras, 
then a comparison of these formulas will imply a relationship between 
mixed volumes and mixed multiplicities. 

\begin{Theorem} \label{relation}
Let $A = k[x_0,...,x_n]$ and $M_0,M_1,...,M_s$ a sequence of sets of 
monomials such that $M_0 = \{x_0,...,x_n\}$ and each $M_i$ consists of 
monomials of the same degree $d_i$ for $i=0,1, \ldots, s.$
Let $\m$ be the maximal graded ideal of $A$ and $J_i$ the ideal generated 
by the monomials of $M_i$. Let $R = R(\m|J_1,...,J_s)$ and let $\Q$ be the 
sequence of polytopes 
$Q_{M_0},Q_{M_1},...,Q_{M_s}$. Then $\deg P_R(u) = n-1$ and
for any $\a \in \NN^{s+1}$ with $|\a| = n-1$,
$$e_\a(\m|J_1,...,J_s) = \frac{MV_{n-1}(\Q_\a)}{\sqrt{n}}.$$
\end{Theorem}

\begin{pf}
Let $S$ denote the subalgebra of the polynomial ring $A[t_0,t_1,...,t_s]$ generated by 
all monomials of the form $f_it_i$ with $f_i \in M_i$. Then $S$ is a 
standard $\NN^{s+1}$-graded algebra over $k$. We shall see that $R \cong 
S$ as $\NN^{s+1}$-graded algebras.
Let $u = (u_0,u_1,..,u_s)$ be any sequence of non-negative integers. The 
vector space 
$R_u$ has a basis  consisting of the monomials of
$\m^{u_0}J_1^{u_1}...J_s^{u_s}$ which are  not contained in  
$\m^{u_0+1}J_1^{u_1}...J_s^{u_s}$.
Since each $J_i$ is generated by $M_i$ and since $M_i$ consists of 
monomials of the same degree, these monomials are of the form 
$f_0f_1...f_s$, where each $f_i$ is a product of $u_i$ 
monomials of $M_i$, $i = 0,1,...,s$. By mapping the elements  
$f_0f_1...f_s \in R_u$ to the elements 
$(f_0t_0^{u_0})(f_1t_1^{u_1})...(f_nt_n^{u_n}) \in S_u$ we obtain an 
$\NN^{s+1}$-graded isomorphism of $R$ and  $S$.\par

Let $\l = (\l_0,\l_1,...,\l_s)$ be any sequence of $s+1$ positive 
integers. 
The above isomorphism induces an $\NN$-graded isomorphism of diagonal 
subalgebras $R^\l \cong S^\l$. 
Let $M^\l$ denote the set of 
all products of $\l_0$ monomials of $M_0$,..., 
$\l_s$ monomials in $M_s$. 
Then
$$S^\l  \cong  k[M^\l].$$
Let $f$ be a product of $s$  monomials 
of $M_1,...,M_s$.  Put $E = 
\{x_1f,...,x_nf\} \subseteq M^{(1,...,1)}$. Then $\ZZ(E)$ contains all 
lattice points of the form $e_i-e_j$, 
where $e_1,...,e_n$ denote the basic vectors of $\RR^n$. 
Therefore, $\ZZ(E)$ contains all lattice points of the hyperplane $x_1 
+ \cdots + x_n = \deg f +1$. Since all monomials of $M^{(1,...,1)}$ have the 
degree $\deg f +1$, the lattice points of $E$ form a basis for 
$\ZZ(M^{(1,...,1)})$. 
Hence $\rank \ZZ(M^{(1,...,1)}) = n$.
Since $Q_E$ is congruent to the convex polytope spanned by the points $e_i$,  
$$V_{n-1}(Q_E) = \frac{\sqrt{n}}{(n-1)!}.$$
Applying Corollary \ref{product} we get
$$e(S^\l) = \frac{(n-1)!}{\sqrt{n}}\sum_{\a \in \NN^{s+1}, |\a| = 
n-1}\frac{1}{\a!}MV_{n-1}(\Q_\a)\l^\a.$$\par

On the other hand, since $ \dim S^\l = \rank \ZZ(M^\l) = n$, 
using Lemma \ref{diagonal} we get $\deg P_S(u) = n-1$ and
$$e(S^\l) = (n-1)!\sum_{\a \in \NN^{s+1}, |\a| = 
n-1}\frac{1}{\a!}e_\a(S)\l^\a.$$
Since the above two formulas for $e(S^\l)$ hold for all sequences $\l$ of
positive integers, we can conclude that their corresponding terms are 
equal. This means 
$$e_\a(S) = \frac{MV_{n-1}(\Q_\a)}{\sqrt{n}}$$
for any $\a \in \NN^{s+1}$ with $|\a| = n-1$.
\end{pf}

It is now easy to interpret mixed volumes as mixed multiplicities of ideals.

\begin{Corollary} \label{main2}
Let $Q_1,...,Q_n$ be an arbitrary collection of lattice convex polytopes in 
$\RR^n$.
Let $A = k[x_0,x_1,...,x_n]$ and let $\m$ be the maximal graded ideal of $A$.
Let $M_i$ be any set of monomials of the same degree in $A$ such that 
$Q_i$ is the convex hull of the lattice points of their dehomogenized 
monomials in $k[x_1,...,x_n]$. Let $J_i$ be the ideal of $A$ generated by 
the monomials  of  $M_i$.  Then
$$MV_n(Q_1,...,Q_n) = e_{(0,1,...,1)}(\m|J_1,...,J_n).$$
\end{Corollary}

\begin{pf}
By definition, the projection of the lattice point of a monomial on the 
hyperplane $x_0 = 0$ is the lattice point of its dehomogenized monomial. 
Therefore,
the convex hull $Q_{M_i}$ of the lattice points of $M_i$ is the projection 
of the polytope $Q_i$ on the hyperplane $x_0 = 0$. As a consequence, the 
volume $V_n(Q_{M_i})$ is proportional to $V_n(Q_i)$. This proportion can be 
computed as
the volume of the convex hull $Q_E$ of the basic vectors 
$e_0,...,e_n$ of $\RR^{n+1}$. Since $V_n(Q_E) = \sqrt{n+1}$, we obtain
$$V_n(Q_i) = \frac{V_n(Q_{M_i})}{\sqrt{n+1}}.$$
From this it follows that the corresponding mixed volumes are also 
proportional:
$$MV_n(Q_1,...,Q_n) = \frac{MV_n(Q_{M_1},...,Q_{M_n})}{\sqrt{n+1}}.$$
On the other hand, applying Theorem \ref{relation} to the sequence 
$M_0,M_1,...,M_n$ of monomials in $n+1$ variables we obtain
$$e_{(0,1,...,1)}(\m|J_1,...,J_n) = \frac{MV_n(Q_{M_1},...,Q_{M_n})}
{\sqrt{n+1}}.$$
Therefore, we can conclude that 
$MV_n(Q_1,...,Q_n) = e_{(0,1,...,1)}(\m|J_1,...,J_n)$.
\end{pf}

An immediate consequence of the interpretation of mixed volumes as 
mixed multiplicities is the non-trivial fact that mixed volumes are 
always non-negative numbers. In fact, we can reprove the following result
given in [Fu2, p. 117].

\begin{Corollary} \label{monoton}
Let $P_1,...,P_n$ and $Q_1,...,Q_n$ be two sequences of convex 
lattice polytopes in $\RR^n$ with $P_i \supseteq Q_i$. Then
$$MV_n(P_1,...,P_n) \ge MV_n(Q_1,...,Q_n).$$
\end{Corollary}

\begin{pf}
By Corollary \ref{main2} we have
\begin{align*}
MV_n(P_1,...,P_n) & = e_{(0,1,...,1)}(\m|I_1,...,I_n),\\
MV_n(Q_1,...,Q_n) & = e_{(0,1,...,1)}(\m|J_1,...,J_n),
\end{align*}
where $I_i$ and $J_i$ are ideals generated by monomial ideals with the 
same degree and $I_i  \supseteq J_i$. Note that the vector space 
$\m^{u_0}I_1^{u_1}...I_n^{u_n}/\m^{u_0+1}I_1^{u_1}...I_n^{u_n}$ 
contains the vector space $\m^{u_0}J_1^{u_1}...J_n^{u_n}/\m^{u_0+1}J_1^{u_1}...J_n^{u_n}$
for all $u = (u_0,u_1,...,u_n) \in \NN^{n+1}$. Then
$$H_{R(\m|I_1,...,I_n)}(u) \ge H_{R(\m|J_1,...,J_n)}(u).$$
Since $e_{(0,1,...,1)}(\m|I_1,...,I_n)$ and 
$e_{(0,1,...,1)}(\m|J_1,...,J_n)$ are the coefficients of one of 
the leading terms  of the corresponding Hilbert polynomials, we obtain 
$$e_{(0,1,...,1)}(\m|I_1,...,I_n) \ge e_{(0,1,...,1)}(\m|J_1,...,J_n),$$
which implies the conclusion.
\end{pf}

\noindent{\bf Remark.} 
{\rm Relations among mixed volumes of lattice polytopes always hold for 
arbitrary convex polytopes by approximating them with  rational convex 
polytopes and then using finer lattices [Te3].}
\smallskip

Now we come to the famous Alexandroff-Fenchel inequality between mixed 
volumes [Fe]:
$$MV_n(Q_1,...,Q_n)^2 \ge MV_n(Q_1,Q_1,Q_3,...,Q_n)MV_n(Q_2,Q_2,Q_3,...,Q_n).$$
Khovanski [Kh] and Teissier [T3] used the Hodge index theorem in 
intersection theory to prove this inequality.
This leads us to believe that a similar inequality should hold between 
mixed multiplicities. \smallskip

\begin{Question} \label{Hodge}
{\rm Let $(A,\m)$ be a local (or standard graded) ring 
with $\dim A = n+1 \ge 3$. Let $I$ be an $\m$-primary ideal and 
$J_1,...,J_n$ ideals of height $n$. Put $\alpha=(0,1, \ldots, 1).$
Is it true that
$$
e_{\alpha}(I|J_1,...,J_n)^2 \ge  
e_{\alpha}(I|J_1,J_1,J_3,...,J_n)e_{\alpha}(I|J_2,J_2,J_3,...,J_n)~?
$$}
\end{Question}

Using Theorem \ref{positivity} we can reduce this theorem to the case 
$\dim A = 3$.
In this case, we have to prove the simpler formula:
$$e_{(0,1,1)}(I|J_1,J_2)^2 \ge e_{(0,1,1)}(I|J_1,J_1)e_{(0,1,1)}(I|J_2,J_2).$$
Unfortunately, we were unable to give an answer to the above question.
The difficulty can be seen from the following observation.\smallskip

\noindent{\bf Remark.}
The above inequality does not hold if $J_1,...,J_n$ are 
$\m$-primary ideals. In this case,  we can even show that the 
inverse inequality holds, namely,
$$e_{\alpha}(I|J_1,...,J_n)^2 \le  
e_{\alpha}(I|J_1,J_1,J_3,...,J_n)e_{\alpha}(I|J_2,J_2,J_3,...,J_n)$$
where $\alpha=(0,1,1, \ldots, 1).$
Using Corollary \ref{Teissier} we can translate it to the inequality
$$e_{(1,1)}(J_1|J_2)^2 \le e(J_1,A)e(J_2,A)$$
for a two-dimensional ring $A$. This inequality was proved first by 
Teissier [Te2] for reduced Cohen-Macaulay rings over an algebraically closed 
field of characteristic zero  and then by Rees and Sharp [RS] in general. 
\smallskip

It is known that computing  mixed volumes is a hard enumerative problem 
(see [EC], [HS1], [HS2] for algorithms and softwares for doing these 
computations). Instead of that we can now compute mixed multiplicities of 
the associated graded ring of the multigraded Rees algebra 
$A[J_1t_1,...,J_nt_n]$ with respect to the ideal $\m$. By 
Corollary \ref{main1}, these mixed multiplicities can be interpreted as 
Samuel multiplicities. The computation of these multiplicities can be 
carried out by computer algebra systems such as {\it Cocoa, Macaulay 2} and 
{\it Singular.} 

\section{Bernstein's theorem}

Let $k[x_1^{\pm 1},...,x_n^{\pm 1}]$ be a Laurent polynomial ring over  
a field $k$.
For any Laurent polynomial 
$$f = \sum_{a \in \ZZ^n}c_ax^a\ (c_a \in k) $$
we will denote by $M(f)$ the set of monomials $x^a$ with $c_a \neq 0$. 
Let 
$Q_f$ denote the convex hull of the lattice points $a$ with $c_a \neq 
0$ in $\RR^n$, i.e.
$Q_f = Q_{M(f)}$. Once calls $Q_f$ the {\it Newton polytope} of $f$. 
\par

Bernstein's theorem  says that  the mixed volume of the associated Newton 
polytopes of $n$ Laurent polynomials is a sharp bound for the number of 
common zeros in the torus $(\CC^*)^{n}$ [Be, Theorem A]. Here we will prove 
Bernstein's theorem by purely algebraic means for  any algebraically 
closed field $k.$

\begin{Theorem} \label{Bernstein}  
Let $k$ be an algebraically closed field. 
Let $f_1,...,f_n$ be Laurent polynomials in $k[x_1^{\pm 1},...,x_n^{\pm 1}]$  
with  finitely many common zeros in $(k^*)^n.$  Then the number of common 
zeros of $f_1,...,f_n$ in $(k^*)^n$ is bounded above by 
$MV_n(Q_{f_1},...,Q_{f_n})$. Moreover, this bound is attained for a 
generic choice of coefficients in  $f_1,..., f_n$  if $k$ has 
characteristic zero.
\end{Theorem}

Here, a generic choice of coefficients in  $f_1,..., f_n$ means that the 
supporting monomials of $f_1,..., f_n$ remain the same while their 
coefficients vary in a non-empty open parameter space. \par

Now we are going to give a homogeneous version of Bernstein's 
theorem.\par

Let $f^h$ denote the homogenization of a Laurent polynomial $f$ in 
$k[x_0^{\pm 1},x_1^{\pm 1},...,x_n^{\pm 1}]$. Then $Q_{f^h}$ is a polytope in 
$\RR^{n+1}$. Its projection 
to the hyperplane $x_0 = 0$ is a polytope  canonically 
identified with $Q_f$. We  have 
$V_n(Q_{f^h}) = \sqrt{n+1}\ V_n(Q_f).$
Hence
$$MV_n(Q_{f_1^h},...,Q_{f_n^h}) =  \sqrt{n+1}\ 
MV_n(Q_{f_1},...,Q_{f_n}).$$
It is also obvious that the number of common zeros of $f_1,...,f_n$ in 
$(k^*)^n$ is equal to the number of common zeros of $f_1^h,...,f_n^h$ 
in $\PP_{k^*}^n$,
where $\PP_{k^*}^n$ denote the set of all points of the projective 
space $\PP_k^n$ with non-zero components.
Thus, Theorem \ref{Bernstein} can be translated as follows.

\begin{Theorem} \label{homogen} 
Let $k$ be an algebraically closed field. Let $g_1,..., g_n$ be homogeneous 
Laurent polynomials in $k[x_0^{\pm 1},x_1^{\pm 1},...,x_n^{\pm 1}]$  with 
finitely many common zeros in $\PP_{k^*}^n.$  Then 
$$ | \{ \a \in \PP_{k^*}^n \mid g_i(\a)=0,\; i=1, 2, \ldots,n \}|
   \leq \frac{MV_n(Q_{g_1},...,Q_{g_n})}{\sqrt{n+1}}.$$
Moreover, this bound is attained for a generic choice of coefficients in  
$g_1,..., g_n$  if $k$ has characteristic zero.
\end{Theorem}

We may reduce the above theorems to the case of polynomials. In fact, 
if we multiply the given Laurent polynomials with an appropriate  
monomial, then we will obtain a new system of polynomials. Obviously, the new 
polynomials in $(k^*)^n$ or in $\PP_{k^*}^n$ have the same common 
zeros. Since their Newton polytopes are translations of the old ones, their 
mixed volumes do not change, too.

Now assume that $g_1,...,g_n$ are homogeneous polynomials in $A  = 
k[x_0,x_1,...,x_n]$. Let $M_i$ be the set of monomials occuring in $g_i$.
Let $\m$ be the maximal graded ideal of $A$ and $J_i$ the 
ideals of $A$ generated by $M_i$. 
Put 
$$R = R(\m|J_1,...,J_n).$$
We know by Theorem \ref{relation} that $\deg P_R(u) = n+1$ and
$$e_{(0,1,...,1)}(R) =  \frac{MV_n(Q_{g_1},..., 
Q_{g_n})}{\sqrt{n+1}}.$$

Therefore, Theorem \ref{homogen} follows from the following result.

\begin{Theorem} \label{multiplicity} 
Let $k$ be an algebraically closed field. Let $g_1,..., g_n$ be homogeneous 
polynomials in $k[x_0,x_1,...,x_n]$  with 
finitely many common zeros in $\PP_{k^*}^n$. Then 
$$ 
   | \{ \a \in \PP_{k^*}^n \mid g_i(\a)=0,\; i=1, 2, \ldots,n \}|
      \leq  e_{(0,1,...,1)}(R). 
$$
Moreover, this bound is attained for a 
generic choice of coefficients in  $g_1,..., g_n$ if $k$ has 
characteristic zero.
\end{Theorem}

\begin{pf}
Let $Q$ be the ideal $(g_1,...,g_n)$. Then there is an 
one-to-one correspondence between common zeros of $g_1,...,g_n$  in 
$\PP_{k^*}^n$ and the one-dimensional homogeneous primes of $A$ which 
contain $Q:(x_0...x_n)^\infty$. As a consequence, the 
assumption on $g_1,...,g_n$ implies that $Q:(x_0...x_n)^\infty$ is a 
one-dimensional ideal. Therefore, the number of common zeros of 
$g_1,...,g_n$  in $\PP_{k^*}^n$ is equal to the number of  minimal 
associated prime ideals of 
$Q:(x_0...x_n)^\infty$ which is bounded above by the multiplicity 
$e(A/(Q:(x_0...x_n)^\infty))$ in view of the associativity formula for 
multiplicities. By the principle of conservation of number 
(see e.g. Fulton [Fu1, Section 10.2]), we only need to show that 
for a generic choice of the coefficients of $g_1,...,g_n,$ 
$Q:(x_0...x_n)^\infty$ is a radical ideal with
$$e(A/(Q:(x_0...x_n)^\infty)) = e_{(0,1,...,1)}(R).$$
Let $J := J_1...J_n$. We  may multiply $g_1,...,g_n$ 
with $x_0...x_n$ to obtain a new system of equations with
$J \subseteq (x_0...x_n)$. 
Since $(x_0...x_n)^m \in J$ for $m \gg 0$, 
$$Q:(x_0...x_n)^\infty = Q:J^\infty.$$
By Corollary \ref{main1} we have
for a generic choice of the coefficients of $g_1,...,g_n,$ 
$$e(A/(Q:J^\infty)) = e_{(0,1,...,1)}(R).$$
Thus, the number of common zeros of $g_1,...,g_n$ in $\PP_{k^*}^n$ is 
bounded by the mixed multiplicity
$e_{(0,1,...,1)}(R)$. 
It remains to show that $Q:J^\infty$ is a radical ideal for a generic 
choice of the coefficients of $g_1,...,g_n$ if $k$ has characteristic zero. 
But this follows from Bertini theorem [Fl, Satz 5.4(e)]. \end{pf}

Finally, we would like to remark that the last statement of the above 
theorems does not hold if the ground field has positive  characteristic.\smallskip

\noindent{\bf Example.} 
Let $k$ be an algebraically closed field with char$(k) = p$. Let 
$f(x) = ax^p + b$ be a polynomial in one variable, $a, b \in k$. 
For $a, b \neq 0$ we choose $c \in k$ such that $c^p = b/a$. Then 
$f(x) = a(x+c)^p$ has only one zero in $k^*$, whereas the Newton polygon 
of $f$ has volume $p$.


\begin{thebibliography}{AAAA}

\bibitem [Ba] {Ba} P. B. Bhattacharya, 
{\it The Hilbert function of two ideals,} 
Proc. Cambridge Phil. Soc. 
{\bf 53 }(1957), 568-575.\par

\bibitem [Be] {Be} D. N. Bernstein, 
{\it The number of roots of a system of  equations} 
(Russian),  Funkcional. Anal. i Prilo\v zen. 
{\bf  9 } (1975), no. 3, 1-4. \par

\bibitem [BF] {BF} T. Bonnesen and W. Fenchel,
{\it  Theorie der konvexen K\"orper,} 
Chelsea, New York, 1971.\par

\bibitem [BH] {BH} W. Bruns and J. Herzog,
{\it  Cohen-Macaulay Rings,} Revised Edition,
Cambridge  Studies in Advanced Mathematics, 39. 
Cambridge University Press, Cambridge, 1998.  \par

\bibitem [CHTV] {CHTV} A. Conca, J. Herzog, N.V. Trung and G. Valla,
{\it Diagonal subalgebras of bigraded algebras and embeddings of 
blow-ups of projective spaces,} 
American Journal of Math. 
{\bf 119 } (1997), 859-901.\par

\bibitem [CLO] {CLO} D. Cox, J. Little and D. O'Shea, 
{\it Using Algebraic Geometry,}
Springer, New York, 1998.\par

\bibitem [EC] {EC} I. Emiris and J. Canny, 
{\it Efficient incremental algorithms for  the sparse resultant and 
the mixed volume,} 
J. Symbolic Comput. 
{\bf 20 } (1995), 117-149.\par

\bibitem [Ew] {Ew} G. Ewald, 
{\it Combinatorial convexity and algebraic geometry,} 
Springer, New York, 1996.\par

\bibitem [Fl] {Fl} H. Flenner,
{\it  Die S\"atze von Bertini f\"ur lokale Ringe,} 
Math. Ann. {\bf 229 } (1977), 97-111. \par

\bibitem [Fu1] {Fu1} W. Fulton, 
{\it Intersection theory,} 
Springer-Verlag,  Berlin-Heidelberg, 1984.\par

\bibitem [Fu2] {Fu2} W. Fulton,
{\it  Introduction to toric varieties,} 
Annals of Mathematics Studies, 131,  
Princeton University Press, 1993.\par

\bibitem [GKZ] {GKZ} I. M. Gelfand, M. Kapranov and A. Zelevinsky,
{\it  Discriminants, Resultants and Multidimensional Determinants,} 
Birkh\"auser, Boston, 1994.\par

\bibitem [HS1] {HS1} B. Huber and B. Sturmfels,
{\it  A polyhedral method for solving sparse polynomial equations,} 
Math. of Computation 
{\bf 64 } (1995), 1541-1555.\par

\bibitem [HS2] {HS2} B. Huber and B. Sturmfels,
{\it  Bernstein's theorem in affine spaces,} 
Discrete. Comput. Geom. 
{\bf 19} (1997), 137-141.\par

\bibitem [KaMV] {KaMV} D. Katz, S. Mandal and J. Verma, 
{\it Hilbert function of bigraded algebras,} 
in:  A. Simis, N. V. Trung and G. Valla  (eds.), 
Commutative Algebra (ICTP, Trieste, 1992), 
291-302, World Scientific, 1994.\par

\bibitem [KaV] {KaV} D. Katz and J. Verma, 
{\it Extended Rees algebras and mixed multiplicities,} 
Math. Z. 
{\bf 202} (1989), 111-128.\par

\bibitem [Kh] {Kh} A. G. Khovanski, 
{\it   Newton polytopes and toric varieties,} 
Functional Anal. Appl. 
{\bf 11} (1977), 289-298.\par

\bibitem [Ku] {Ku} A. G. Kuschnirenko, 
{\it Newton polytopes and Bezout theorem,} 
Functional Anal. Appl. 
{\bf 10} (1976), 233-235.\par

\bibitem [R1] {R1} D. Rees,
{\it  $\frak a$-transforms of local rings and a theorem on  
multiplicities of ideals,} 
Proc. Cambridge Philos. Soc. 
{\bf 57} (1961), 8-17. \par

\bibitem [R2] {R2} D. Rees, 
{\it Generalizations of reductions and mixed multiplicities,} 
J. London Math. Soc. 
{\bf 29} (1984), 423-432. \par

\bibitem [Ro] {Ro} P. Roberts, 
{\it Local Chern classes, multiplicities and perfect complexes,} 
M\'emoire Soc. Math. France 
{\bf 38} (1989), 145-161. \par

\bibitem [RS] {RS} D. Rees and R. Y. Sharp,
{\it  On a theorem of B. Teissier on mixed  multiplicities of 
ideals in local rings,} 
J. London Math. Soc. 
{\bf 18} (1978), 449-463.\par

\bibitem [Sta] {Sta} R. P.  Stanley, 
{\it Combinatorics and Commutative Algebra,} 
Birh\"auser, Boston, 1986. \par

\bibitem [Stu] {Stu} B. Sturmfels, 
{\it Solving systems of polynomial equations,}
CBMS Regional Conference Series in Mathematics, No. 97, 
American Mathematical Society, 2002 .  \par

\bibitem [SV] {SV} J. St\"uckrad and W. Vogel, 
{\it Buchsbaum rings and applications,} 
VEB Deutscher Verlag der Wisssenschaften, Berlin, 1986.\par

\bibitem [Sw] {Sw} I. Swanson, {\it Mixed multiplicities, joint reductions and quasi-unmixed local rings,} J. London Math. Soc. {\bf 48} (1993), 1-14.\par 

\bibitem [Te1] {Te1} B. Teissier, 
{\it  Cycles \'evanescents, sections planes, et 
conditions de Whitney,} 
Singularit\'es \`a Carg\`ese 1972, Ast\`erisque 
{\bf 7-8 }(1973), 285-362.  \par

\bibitem [Te2] {Te2} B. Teissier, 
{\it Sur un in\'egalit\'e \`a la Minkowski pour les  multiplicit\'es,} 
(Appendix to a paper by D. Eisenbud and H. I. Levine),  
Ann. of Math. 
{\bf 106} (1977), 38-44.\par

\bibitem [Te3] {Te3} B. Teissier, 
{\it Du th\'eor\`eme de l'index de Hodge aux in\'egalit\'es 
isop\'erim\'etriques, } 
C. R. Acad. Sci. Paris Ser. A-B  
{\bf 288 } (1979), no. 4, A287--A289.  \par

\bibitem   [Tr1] {Tr1} N. V. Trung, 
{\it The Castelnuovo regularity of the Rees 
algebra and the associated graded ring,} 
Trans. Amer. Math. Soc. 
{\bf 350} (1998), 2813-2832. \par

\bibitem   [Tr2] {Tr2} N. V. Trung,
{\it  Positivity of mixed multiplicities,}  Math. Ann.  
{\bf 319 } (2001), 33--63.  \par 

\bibitem  [Va] {Va} G. Valla, 
{\it Certain graded algebras are always Cohen-Macaulay,} 
J. Algebra, 
{\bf 42 }(1976), 537-548.\par

\bibitem [Vi] {Vi} D. Q. Viet, {\it Mixed multiplicities of arbitrary ideals in local rings,}
Comm. Algebra {\bf 28} (8) (2000), 3803-3821.\par

\bibitem [Ve1] {Ve1} J. K. Verma, 
{\it Rees algebras and mixed multiplicities,} 
Proc. Amer. Math. Soc. 
{\bf 104} (1988), 1036-1044. \par

\bibitem [Ve2] {Ve2} J. K. Verma, 
{\it Multigraded Rees algebras and mixed multiplicities,} 
J. Pure Appl. Algebra 
{\bf 77 }(1992), 219-228. \par

\bibitem  [Wa] {Wa} B. L. Van der Waerden,
{\it  On Hilbert's function, series of composition of ideals 
and a generalization of the theorem of Bezout,} 
Proc. K. Akad. Wet. Amsterdam 
{\bf 31 }(1928), 749-770.

\end{thebibliography}
\end{document}